\documentclass[runningheads]{llncs}
\usepackage{graphicx}
\usepackage{hyperref}
\usepackage{amsfonts, amsfonts, amssymb, amsmath}
\usepackage{color}
\hypersetup{
    colorlinks=true,
    linkcolor=blue,
    filecolor=magenta,      
    urlcolor=blue,}

\newcommand{\R}{\mathbb{R}}
\newcommand{\cS}{\mathcal S}
\newcommand{\Ver}{\mathrm{Ver}}
\newcommand{\ver}{\mathrm{ver}}
\newcommand{\Hor}{\mathrm{Hor}}
\newcommand{\hor}{\mathrm{hor}}
\newcommand{\Skew}{\mathrm{Skew}}
\newcommand{\Sym}{\mathrm{Sym}}
\newcommand{\Exp}{\mathrm{Exp}}
\newcommand{\Log}{\mathrm{Log}}
\newcommand{\Tr}{\mathrm{Tr}}

\begin{document}
\title{Parallel Transport on Kendall Shape Spaces}
\author{Nicolas Guigui\inst{1}\orcidID{0000-0002-7901-0732} \and
Elodie Maignant\inst{1,2} \and Alain Trouvé \inst{2}
Xavier Pennec\inst{1}}
\authorrunning{N. Guigui et al.}
\institute{Université Côte d'Azur, Inria Epione project team, France \and
Centre Borelli, ENS Paris Saclay}
\maketitle              % typeset the header of the contribution
\begin{abstract}
Kendall shape spaces are a widely used framework for the statistical analysis of shape data arising from many domains, often requiring the parallel transport as a tool to normalise time series data or transport gradient in optimisation procedures. We present an implementation of the pole ladder, an algorithm to compute parallel transport based on geodesic parallelograms and compare it to methods by integration of the parallel transport ordinary differential equation.

\end{abstract}
\section{Introduction}
Kendall shape spaces are a ubiquitous framework for the statistical analysis of data arising from medical imaging, computer vision, biology, chemistry and many more domains. The underlying idea is that a shape is what is left after removing the effects of rotation, translation and re-scaling, and to define a metric accounting for those invariances. This involves defining a Riemannian submersion and the associated quotient structure, resulting in non trivial differential geometries with singularities and curvature, requiring specific statistical tools to deal with such data.

In this context parallel transport is a fundamental tool to define statistical models and optimisation procedures, such as the geodesic or spline regression \cite{kim_smoothing_2020,nava-yazdani_geodesic_2020} and for the normalisation of time series of shapes \cite{lorenzi_efficient_2014,cury_spatio-temporal_2016}. 
However parallel transport is defined by an ordinary differential equation (ODE) and there is usually no closed-form solution.
Approximation methods have therefore been derived, either by direct integration \cite{kim_smoothing_2020}, or by integration of the geodesic equation to approximate Jacobi fields (the fanning scheme \cite{louis_fanning_2018}). Another class of approximations refereed to as ladder methods, relies on iterative constructions of geodesic parallelograms that only require approximate geodesics \cite{guigui_numerical_2020}.

In this work, we present an implementation of the pole ladder that leverages the quotient structure of Kendall shape spaces and strongly relies on the open-source Python package \url{geomstats}. We compare it to the method of Kim et al.\ \cite{kim_smoothing_2020} by approximate integration.

We first recall the quotient structure of Kendall shape spaces and its use by Kim et al. to compute parallel transport in sec.~\ref{sec:kendall}, then in sec.~\ref{sec:pole} we recall the pole ladder scheme and the main result from \cite{guigui_numerical_2020} on its convergence properties. Numerical simulations to compare the two methods are reported in sec.~\ref{sec:results}.

\section{The quotient structure of Kendall shape spaces}
\label{sec:kendall}
We first describe the shape space, as quotient of the space of configurations of $k$ points in $\R^m$ (called landmarks) by the groups of translations, re-scaling and rotations of $\R^m$. For a thorough treatment this topic, we refer the reader to \cite{kendall_shape_1984,dryden_statistical_2016}.

\subsection{The pre-shape space}

We define the space of $k$ landmarks of $\R^m$ as the space of $m \times k$ matrices $M(m,k)$. For $x \in M(m,k)$, let $x_i$ denote the columns of $x$, i.e. points of $\R^m$ and let $\bar x$ be their barycentre. We remove the effects of translation by considering the matrix with columns $x_i - \bar x$ instead of $x$. We further remove the effects of scaling by dividing $x$ by its Frobenius norm (written $\|\cdot\|$). This defines the pre-shape space $\cS_m^k = \{x \in M(m,k)\;|\; \sum_{i=1}^k x_i = 0, \;\; \|x\|=1\}$, which is identified with the hypersphere of dimension $m(k-1) -1$. The pre-shape space is therefore a differential manifold whose tangent space at any $x \in \cS^k_m$ is given by $T_x\cS^k_m=\{w \in M(m, k) \;|\; \sum_{i=1}^k w_i=0,\; \Tr(w^Tx)=0\}$.

The ambient Frobenius metric $\langle\cdot,\cdot\rangle$ thus defines a Riemannian metric on the pre-shape space, with constant sectional curvature and known geodesics: let $x,y \in \cS^k_m$ with $x \neq y$, and $w \in T_x \cS^k_m$
\begin{align}
\label{eq:top_exp}
    x_w=\exp_x(w) &= \cos(\|w\|)x + \sin(\|w\|)\frac{w}{\|w\|},\\
    \log_x(y) &= \arccos(\langle y,x\rangle )\frac{y - \langle y,x\rangle x}{\|y - \langle y,x\rangle x\|}.
    % \\
    % \Pi_x^{x_w}v &= \langle v,\frac{w}{\|w\|}\rangle  \left(- \sin(\|w\|)x + \cos(\|w\|)\frac{w}{\|w\|}\right) \\
    %     &\quad+ \left(v - \langle v,\frac{w}{\|w\|}\rangle \frac{w}{\|w\|} \right).
\end{align}

Moreover, this metric is invariant to the action of the rotation group $SO(m)$. This allows to define the shape space as the quotient $\Sigma_m^k~= S_m^k/SO(m)$.

\subsection{The shape space}
To remove the effect of rotations, we define the equivalence relation $\sim$ on $\cS^k_m$ by $x\sim y \iff \exists R\in SO(m)$ such that $y = Rx$. For $x\in \cS^k_m$, let $[x]$ denote its equivalence class for $\sim$. This equivalence relation results from the group action of $SO(m)$ on $\R^m$. 
This action is smooth, proper but not free everywhere when $m\geq 3$. This makes the orbit space $\Sigma_m^k = \{[x] \;|\; x\in \cS^k_m\}$ a differential manifold with singularities where the action is not free. 

One can describe these singularities explicitly : they correspond to the matrices of $\cS^k_m$ of rank $m-2$ or less \cite{le_riemannian_1993}. For $k\geq 3$, the spaces $\Sigma^k_1$ and $\Sigma^k_2$ are always smooth.
Moreover, as soon as $m\geq k$, the manifold acquires boundaries. As an example, while the space $\Sigma^3_2$ of 2D triangles is identified with the sphere $S^2(1/2)$, the space $\Sigma^3_3$ of 3D triangles is isometric to a 2-ball \cite{le_riemannian_1993}.

Away from the singularities, the canonical projection map $\pi:x \mapsto [x]$ is a Riemannian submersion, and plays a major role in defining the metric on the shape space. Let $d_x\pi$ be its differential map at $x\in \cS^k_m$, whose kernel defines the vertical tangent space, which corresponds to the tangent space of the submanifold $\pi^{-1}([x])$, called fiber above $[x]$:
\begin{equation*}
    \Ver_x = \{Ax \;|\; A \in \Skew(m)\} = \Skew(m)\cdot x
\end{equation*}
where $\Skew(m)$ is the space of skew-symmetric matrices of size $m$.

\subsection{The quotient metric}
The Frobenius metric on the pre-shape space allows to define the horizontal spaces as the orthogonal complements to the vertical spaces:
\begin{align*}
    \Hor_x &= \{w \in T_x\cS^k_m \;|\; Tr(Axw^T)=0 \; \forall A \in \Skew(m)\} \\
       &= \{w \in T_x\cS^k_m \;|\; xw^T \in \Sym(m)\}
\end{align*}
where $\Sym(m)$ is the space of symmetric matrices of size $m$. 
Lemma 1 from \cite{nava-yazdani_geodesic_2020} allows to compute the vertical component of any tangent vector:

\begin{lemma}
For any $x\in \cS^k_m$ and $w \in T_x\cS^k_m$, the vertical component of $w$ can be computed as $\Ver_x(w) = Ax$ where $A$ solves the Sylvester equation:
\begin{equation}
    \label{eq:sylvester}
    Axx^T + xx^TA = wx^T - xw^T
\end{equation}
If $\mathrm{rank}(x) \geq m-1$, $A$ is the unique skew-symmetric solution of \eqref{eq:sylvester}.
\end{lemma}
In practice, the Sylvester equation can be solved by an eigenvalue decomposition of $xx^T$. This defines $\ver_x$, the orthogonal projection on $\Ver_x$.
As $T_x\cS^k_m = \Ver_x \oplus \Hor_x$, any tangent vector $w$ at $x \in \cS^k_m$ may be decomposed into a horizontal and a vertical component, by solving \eqref{eq:sylvester} to compute $\ver_x(w)$, and then $\hor_x(w)=w - \ver_x(w)$.

Furthermore, as $\Ver_x = \ker(d_x\pi)$, $d_x\pi$ is a linear isomorphism from $\Hor_x$ to $T_{[x]}\Sigma^k_m$. The metric on $\Sigma^k_m$ is defined such that this isomorphism is an isometry. Note that the metric does not depend on the choice of the $y$ in the fiber $\pi^{-1}([x])$ since all $y$ in $\pi^{-1}([x])$ may be obtained by a rotation of $x$, and the Frobenius metric is invariant to the action of rotations.
This makes $\pi$ a Riemannian submersion. Additionally, $\pi$ is surjective so for every vector field on $\Sigma^k_m$ there is a unique horizontal lift, i.e.\ a vector field on $\cS^k_m$ whose vertical component is null everywhere. The tangent vectors of $\Sigma^k_m$ can therefore be identified with horizontal vectors of $\cS^k_m$. 
One of the main characteristics of Riemannian submersions was proved by O'Neill \cite{oneill_semi-riemannian_1983}:

\begin{theorem}[O'Neill]
\label{thm}
Let $\pi:M \rightarrow B$ be a Riemannian submersion. If $\gamma$ is a geodesic in $M$ such that $\dot \gamma(0)$ is a horizontal vector, then $\dot \gamma$ is horizontal everywhere and $\pi \circ \gamma$ is a geodesic of $B$ of the same length as $\gamma$.
\end{theorem}

\begin{remark}
\label{main_remark}
We emphasise that an equivalent proposition cannot be derived for the parallel transport of a tangent vector. Indeed the parallel transport of a horizontal vector field along a horizontal geodesic may not be horizontal. This will be detailed in the next subsection and constitutes a good example of metric for which computing geodesics is easier than computing parallel transport, although the former is a variational problem and the latter is a linear ODE.
\end{remark}

Furthermore, the Riemannian distances $d$ on $\cS^k_m$ and $d_\Sigma$ on $\Sigma^k_m$ are related by
\begin{equation}
    d_\Sigma(\pi(x), \pi(y)) = \inf_{R \in SO(m)}d(x, Ry).
\end{equation}
The optimal rotation $R$ between any $x,y$ is unique in a subset $U$ of $\cS^k_m \times \cS^k_m$, which allows to define the \textit{align} map $\omega:U\rightarrow \cS^k_m$ that maps $(x,y)$ to $Ry$. In this case, $ d_\Sigma(\pi(x), \pi(y)) = d(x, \omega(x,y))$ and $x\omega(x,y)^T \in \Sym(m)$. It is useful to notice that $w(x,y)$ can be directly computed by a pseudo-singular value decomposition of $xy^T$\cite{kendall2009statistical}. Finally, $x$ and $\omega(x,y)$ are joined by a horizontal geodesic.

\subsection{Implementation in \url{geomstats}}
The \url{geomstats} library \cite{miolane_geomstats_2020}, available at \url{https://geomstats.ai}, implements classes of manifolds equipped with Riemannian metrics. 
It contains an abstract class for quotient metrics, that allows to compute the Riemannian distance, exponential and logarithm maps in the quotient space from the ones in the top space.

In the case of the Kendall shape spaces, the quotient space cannot be seen as a submanifold of some $\R^N$. Moreover, the projection $\pi$ and its total derivative $d\pi$ can't be computed explicitly.
However, the align map amounts to identifying the shape space with a local horizontal section of the pre-shape space, and thanks to the characteristics of Riemannian submersions mentioned in the previous subsections, all the computations can be done in the pre-shape space.

Recall that $\exp$, $\log$, and $d$ denote the operations of the pre-shape space $\cS^k_m$ and are given in \eqref{eq:top_exp}. We obtain from theorem \ref{thm} for any $x,y \in \cS^k_m$ and $v \in T_x\cS^k_m$
\begin{align*}
    \exp_{\Sigma,[x]}(d_x\pi v) &= \pi(\exp_x(\hor_x(v))), \\
    \log_{\Sigma,[x]}([y]) &= d_x\pi \log_x(\omega(x, y)), \\
    d_{\Sigma}([x], [y]) &= d(x, \omega(x,y)).
\end{align*}

% \subsection{Exponential and logarithm maps, horizontal geodesics.}
% Explain quickly how the exponential and log maps project on lower space, and how one can describe the lower geodesics.
% Explain that this does not work for parallel transport. All this subsection should lead to next section, and account for it.

\subsection{Parallel transport in the shape space}
As noticed in Remark~\ref{main_remark}, one cannot use the projection of the parallel transport in the pre-shape space $\cS^k_m$ to compute the parallel transport in the shape space $\Sigma^k_m$. Indeed \cite{kim_smoothing_2020} proved the following

\begin{proposition}[Kim et al. \cite{kim_smoothing_2020}]
Let $\gamma$ be a horizontal $C^1$-curve in $\cS^k_m$ and $v$ be a horizontal tangent vector at $\gamma(0)$. Assume that $\mathrm{rank}(\gamma(s) \geq m-1$ except for finitely many $s$. Then the vector field $s\mapsto v(s)$ along $\gamma$ is horizontal and the projection of $v(s)$ to $T_{[\gamma(s)]}\Sigma^k_m$ is the parallel transport of $d_x\pi v$ along $[\gamma(s)]$ if and only if $s\mapsto v(s)$ is the solution of
\begin{equation}
    \label{eq:pt}
    \dot v(s) = - \Tr(\dot \gamma(s) v(s)^T) \gamma(s) + A(s)\gamma(s), \;\qquad v(0) = v
\end{equation}
where for every $s$, $A(s)\in \Skew(m)$ is the unique solution to
\begin{equation}
    \label{eq:sylvester_dot}
    A(s)\gamma(s)\gamma(s)^T + \gamma(s)\gamma(s)^TA(s) = \dot \gamma(s)v(s)^T - v(s)\dot\gamma(s)^T.
\end{equation}
\end{proposition}
Eq.~\eqref{eq:pt} means that the covariant derivative of $s\mapsto v(s)$ along $\gamma$ must be a vertical vector at all times, defined by the matrix $A(s) \in \Skew(m)$. These equations can be used to compute parallel transport in the shape space. To compute the parallel transport of $d_x\pi w$ along $[\gamma]$, \cite{kim_smoothing_2020} propose the following method: one first chooses a discretization time-step $\delta = \frac{1}{n}$, then repeat for every $s=\frac{i}{n}, i=0\ldots n$
\begin{enumerate}
    \item Compute $\gamma(s)$ and $\dot \gamma(s)$,
    \item Solve the Sylvester equation \eqref{eq:sylvester_dot} to compute $A(s)$ and the r.h.s. of $\eqref{eq:pt}$,
    \item Take a discrete Euler step to obtain $\tilde v(s + \delta)$
    \item Project $\tilde v(s + \delta)$ to $T_{\gamma(s)}\cS^k_m$ to obtain $\hat v(s + \delta)$,
    \item Project to the horizontal subspace: $v(s+\delta) \leftarrow \hor(\hat v(s + \delta))$
    \item $s\leftarrow s + \delta$
\end{enumerate}
We notice that this method can be accelerated by a higher-order integration scheme, such as Runge-Kutta (RK) by directly integrating the system $\dot v = f(v, s)$ where $f$ is a smooth map given by \eqref{eq:pt} and \eqref{eq:sylvester_dot}. In this case, steps 4. and 5. are not necessary.
The precision and complexity of this method is then bound to that of the integration scheme used. As ladder methods rely only on geodesics, which can be computed in closed-form and their convergence properties are well understood \cite{guigui_numerical_2020}, we compare this method by integration to the pole ladder. We focus on the case where $\gamma$ is a horizontal geodesic.

\section{The Pole ladder algorithm}
\label{sec:pole}

\subsection{Description}
The pole ladder is a modification of the Schild's ladder \cite{misner_gravitation_1973} proposed by \cite{lorenzi_efficient_2014}. The pole ladder is more precise and cheaper to compute as shown by \cite{guigui_numerical_2020}. It is also exact in symmetric spaces \cite{pennec_parallel_2018}. We thus focus on this method. We describe it here in a Riemannian manifold $(M, \langle,\rangle )$.

Consider a geodesic curve $\gamma: t\mapsto \gamma(t) \in M$, with initial conditions $x~=\gamma(0) \in M$ and $w = \dot\gamma(0) \in T_xM$. In order to compute the parallel transport of $v \in T_xM$ along $\gamma$, between times $x$ and $y=\gamma(1)$, the pole ladder consists in first dividing the main geodesic $\gamma$ in $n$ segments of equal length and computing the geodesic from $x$ with initial velocity $\frac{v}{n^\alpha}$, obtaining $x_v = \Exp_x(\frac{v}{n^\alpha})$. Then for each segment to repeat the following construction (see figure~\ref{fig:pole_construction}):

\begin{figure}[t]
    \centering
    \includegraphics[width=0.5\textwidth]{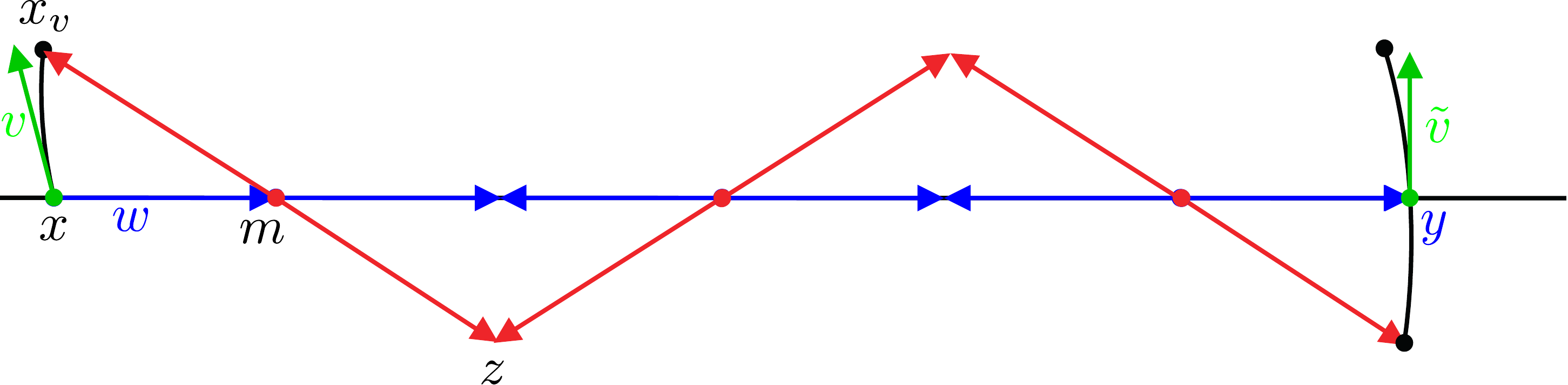}
    \vspace{-4mm}
    \caption{Schematic representation of the pole ladder}
    \label{fig:pole_construction}
    \vspace{-4mm}
\end{figure}

\begin{figure}[b]
    \begin{minipage}{.3\textwidth}
        \centering
        \includegraphics[width=0.8\textwidth]{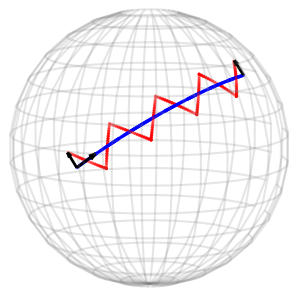}
    \end{minipage}
    \begin{minipage}{.33\textwidth}
        \centering
        \includegraphics[width=0.99\textwidth]{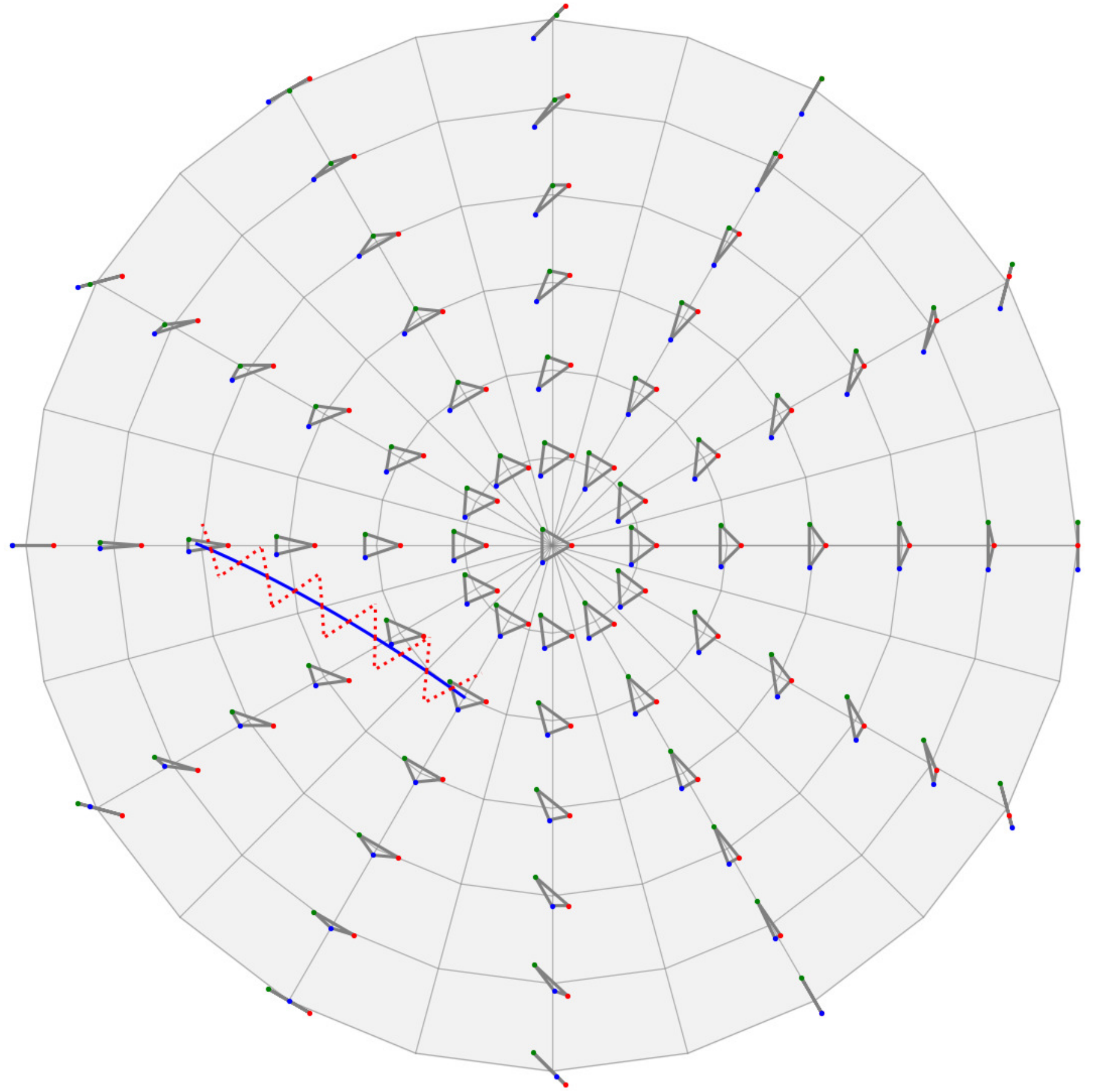}
    \end{minipage}
    \begin{minipage}{.33\textwidth}
        \centering
        \includegraphics[width=0.99\textwidth]{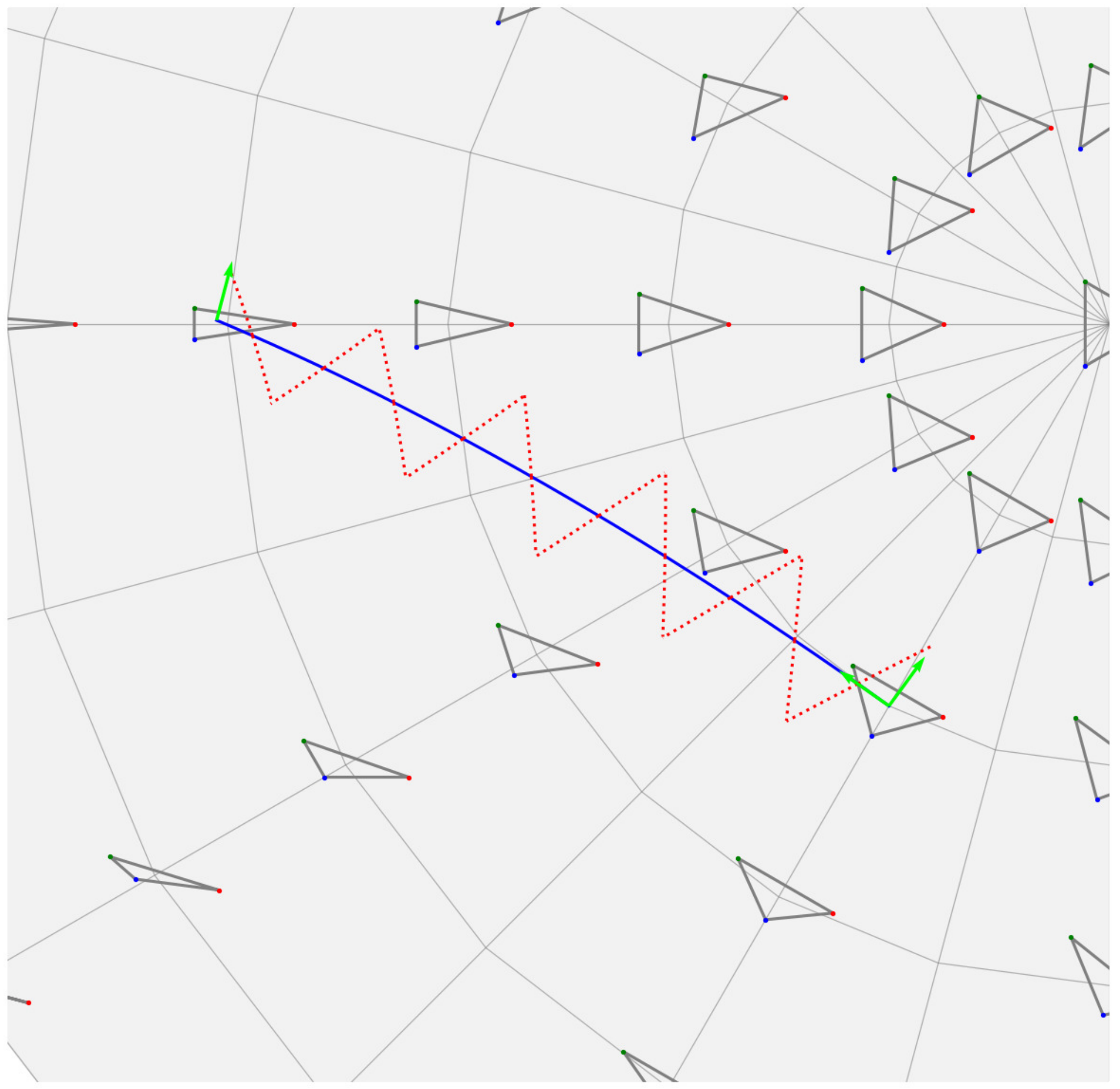}
    \end{minipage}
    \caption{Visualisation of the pole ladder on $S^2$ (left) and $\Sigma_3^3$ (middle and right)}
    \label{fig:pole_visu_kendall}
\end{figure}

\begin{enumerate}
    \item Compute the midpoint of the segment $m = \Exp_x(\frac{w}{2n})$ and the initial speed of the geodesic from $m$ to $x_v$: $a=\Log_m(x_v)$.
    \item Extend this diagonal geodesic by the same length to obtain $z=\Exp_m(-a)$.
    \item Repeat steps 2 and 3 with $x_v \leftarrow z$ and $m \leftarrow \Exp_m(\frac{w}{n})$.
\end{enumerate}
After $n$ steps, compute $\tilde v = n^\alpha(-1)^n \Log_y(z)$. According to \cite{guigui_numerical_2020}, $\alpha\geq 1$ can be chosen, and $\alpha=2$ is optimal. This vector is an approximation of the parallel transport of $v$ along $\gamma$, $\Pi_x^{x_w}v$. This is illustrated on the 2-sphere and in the case $k=m=3$ on Figure~\ref{fig:pole_visu_kendall}.

% \begin{figure}[b]
%     \centering
%     \includegraphics[width=0.35\textwidth]{pole_visu.eps}
%     \caption{Visualisation of the pole ladder on the $2$-sphere}
%     \label{fig:pole_visu}
% \end{figure}

\subsection{Properties}
The pole ladder is studied in depth in \cite{guigui_numerical_2020}. We give here the two main properties. Beside its quadratic convergence speed, the main advantage is that this method is available as soon as geodesics are known (even approximately). It is thus applicable very easily in the case of quotient metrics.
\begin{theorem}
\begin{itemize}
    \item The pole ladder converges to the exact parallel transport when the number of steps $n$ goes to infinity, and the error decreases in $O(\frac{1}{n^2})$, with rate related to the covariant derivative of the curvature tensor.
    \item If $M$ is a symmetric space, then the pole ladder is exact in just one step.
\end{itemize}
\end{theorem}
For instance, $\Sigma^3_2$ is symmetric, making pole ladder exact in this case.

\subsection{Complexity}
The main drawback of ladder schemes is that logarithms are required. Indeed the Riemannian logarithm is only locally defined, and often solved by an optimisation problem when geodesics are not known in closed form. 

In the case of Kendall shape spaces, it only requires to compute an alignment step, through a singular value decomposition, with usual complexity $O(m^3)$, then the $\log$ of the hypersphere, with linear complexity.
Moreover, the result of $\log\circ \;\omega$ is horizontal, so the vertical component needs not be computed for the exponential of step 2, and only the $\exp$ of the hypersphere, also with linear complexity, needs to be computed.
The vertical projection needs to be computed for the first step. Solving the Sylvester equation through an eigenvalue decomposition also has complexity $m^3$.
For $n$ rungs of the pole ladder, the overall complexity is thus $O((n+1)(m^3 + 2mk)) + mk + m^3) = O(nm^3)$.

On the other hand, the method by integration doesn't require logarithms but requires solving a Sylvester equation and a vertical decomposition at every step. The overall complexity is thus $O(2nm^3 + mk)$. Both algorithms are thus comparable in terms of computational cost for a single step.

\section{Numerical Simulations and Results}
\label{sec:results}
\begin{figure}[t]
    \begin{minipage}{.49\textwidth}
        \centering
        \includegraphics[width=0.99\textwidth]{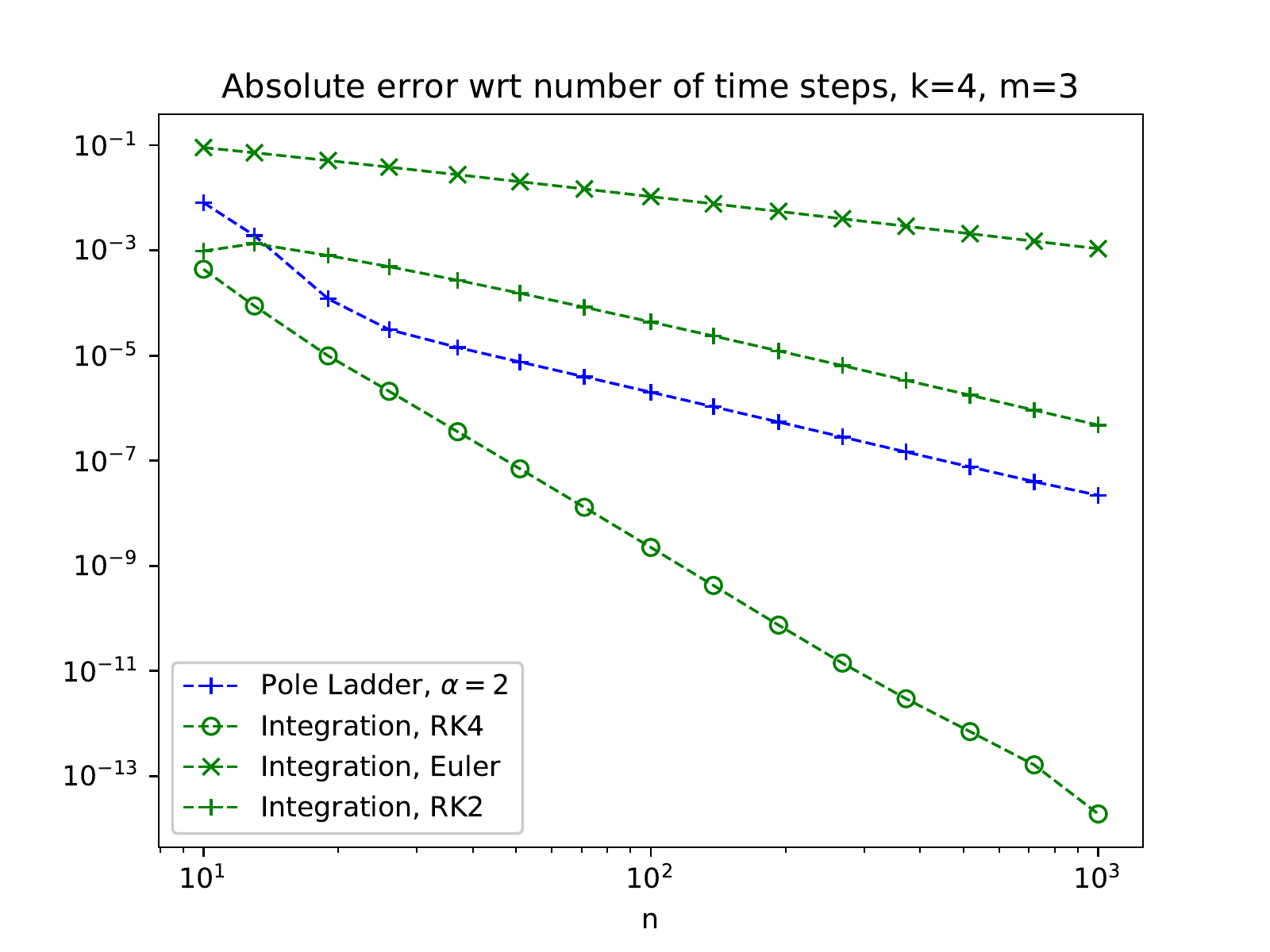}
    \end{minipage}
    \begin{minipage}{.5\textwidth}
        \centering
        \includegraphics[width=0.99\textwidth]{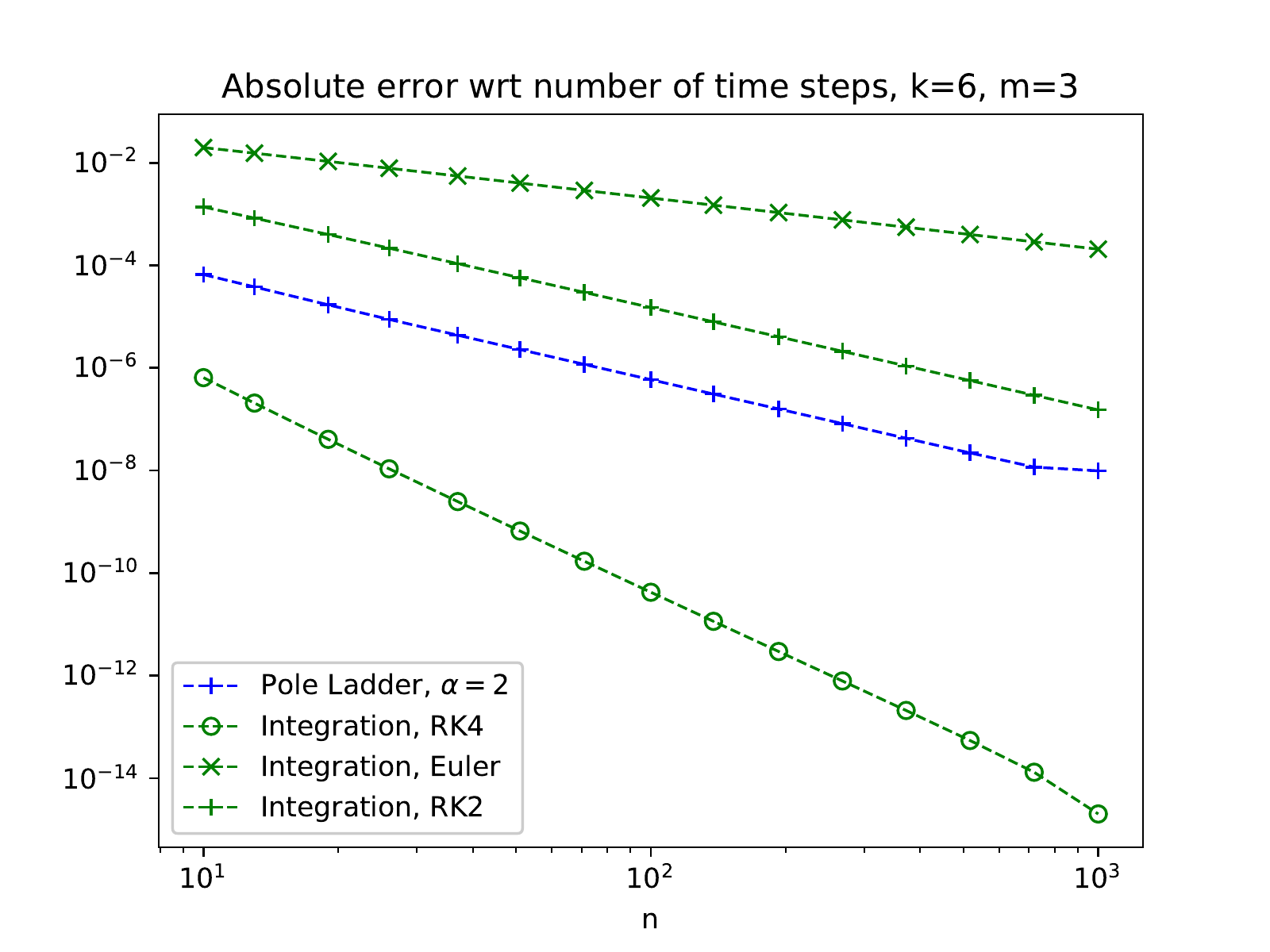}
    \end{minipage}
    \vspace{-6mm}
    \caption{Error of the parallel transport of $v$ along the geodesic with initial velocity $w$ where $v$ and $w$ are orthonormal.}
     \vspace{-4mm}
   \label{fig:result}
\end{figure}

We draw a point $x$ at random in the pre-shape space, along with two orthogonal horizontal unit tangent vectors $v,w$, and compute the parallel transport of $d_x\pi v$ along the geodesic with initial velocity $d_x\pi w$. We use a number of steps $n$ between $10$ and $1000$ and the result with $n=1100$ as the reference value to compute the error made by lower numbers of steps. The results are displayed on Figure~\ref{fig:result} for the cases $k=4,6$ and $m=3$ in log-log plots.
As expected, the method proposed by \cite{kim_smoothing_2020} converges linearly, while RK schemes of order two and four show significant acceleration. 
The pole ladder converges with quadratic speed and thus compares with the RK method of order two, although the complexity of the RK method is multiplied by its order.

\section{Conclusion and future work}
We presented the Kendall shape space and metric, highlighting the properties stemming from its quotient structure. 
This allows to compute parallel transport with the pole ladder using closed-form solution for the geodesics. 
% The comparison with a method by integration showed that pole ladder reached a given precision with a much lower number steps, for a comparable computational cost.
This off-the-shelf algorithm can now be used in learning algorithms such as geodesic regression or local non-linear embedding. This will be developed in future works.

\section{Acknowledgments}
\label{sec:acknowledgments}
This work was partially funded by the ERC grant Nr. 786854 G-Statistics from the European Research Council under the European Union’s Horizon 2020 research and innovation program.
It was also supported by the French government through the 3IA Côte d’Azur Investments ANR-19-P3IA-0002 managed by the National Research Agency.

\bibliographystyle{splncs04}
%\bibliography{Kendall}

\begin{thebibliography}{}
\providecommand{\url}[1]{\texttt{#1}}
\providecommand{\urlprefix}{URL }
\providecommand{\doi}[1]{https://doi.org/#1}

\end{thebibliography}


\begin{thebibliography}{10}
\providecommand{\url}[1]{\texttt{#1}}
\providecommand{\urlprefix}{URL }
\providecommand{\doi}[1]{https://doi.org/#1}

\bibitem{cury_spatio-temporal_2016}
Cury, C., Lorenzi, M., Cash, D., Nicholas, J., Routier, A., Rohrer, J.,
  Ourselin, S., Durrleman, S., Modat, M.: Spatio-{Temporal} {Shape} {Analysis}
  of {Cross}-{Sectional} {Data} for {Detection} of {Early} {Changes} in
  {Neurodegenerative} {Disease}. In: {SeSAMI} 2016 - {First} {International}
  {Workshop} {Spectral} and {Shape} {Analysis} in {Medical} {Imaging}. LNCS
  10126, pp. 63 -- 75. Springer (Sep 2016). %\doi{10.1007/978-3-319-51237-2\_6}

\bibitem{dryden_statistical_2016}
Dryden, I.L., Mardia, K.V.: Statistical {Shape} {Analysis}: {With}
  {Applications} in {R}. John Wiley \& Sons (Sep 2016)

\bibitem{guigui_numerical_2020}
Guigui, N., Pennec, X.: Numerical {Accuracy} of {Ladder} {Schemes} for
  {Parallel} {Transport} on {Manifolds} (Jul 2020).
  % \url{https://hal.inria.fr/hal-02894783}

\bibitem{kendall_shape_1984}
Kendall, D.G.: Shape {Manifolds}, {Procrustean} {Metrics}, and {Complex}
  {Projective} {Spaces}. Bulletin of the London Mathematical Society
  \textbf{16}(2),  81--121 (1984).
    % \doi{https://doi.org/10.1112/blms/16.2.81},
  % \url{https://londmathsoc.onlinelibrary.wiley.com/doi/abs/10.1112/blms/16.2.81},
  % \_eprint: https://londmathsoc.onlinelibrary.wiley.com/doi/pdf/10.1112/blms/16.2.81
  
\bibitem{kendall2009statistical}
Kendall, W. S., Le, H.: Statistical {Shape} Theory. New Perspectives in
    Stochastic Geometry, 348-373 (2009).


\bibitem{kim_smoothing_2020}
Kim, K.R., Dryden, I.L., Le, H., Severn, K.E.: Smoothing splines on
  {Riemannian} manifolds, with applications to {3D} shape space. Journal of the
  Royal Statistical Society: Series B (Statistical Methodology) (Dec 2020).
  % \doi{10.1111/rssb.12402}

\bibitem{le_riemannian_1993}
Le, H., Kendall, D.G.: The {Riemannian} {Structure} of {Euclidean} {Shape}
  {Spaces}: {A} {Novel} {Environment} for {Statistics}. The Annals of
  Statistics  \textbf{21}(3),  1225--1271 (1993),
  % \url{https://www.jstor.org/stable/2242196}, publisher: Institute of Mathematical Statistics

\bibitem{lorenzi_efficient_2014}
Lorenzi, M., Pennec, X.: Efficient {Parallel} {Transport} of {Deformations} in
  {Time} {Series} of {Images}: {From} {Schild} to {Pole} {Ladder}. J Mathhttps://www.overleaf.com/project/601fe6b8ac2c5d67d48ca497
  Imaging Vis  \textbf{50}(1),  5--17 (Sep 2014).
  %\doi{10.1007/s10851-013-0470-3}

\bibitem{louis_fanning_2018}
Louis, M., Charlier, B., Jusselin, P., Pal, S., Durrleman, S.: A {Fanning}
  {Scheme} for the {Parallel} {Transport} {Along} {Geodesics} on {Riemannian}
  {Manifolds}. SIAM Journal on Numerical Analysis  \textbf{56}(4), 2563--2584
  (2018). %\doi{10.1137/17M1130617}

\bibitem{miolane_geomstats_2020}
  Miolane, N., Guigui, N., Brigant, A.L., Mathe, J., Hou, B., Thanwerdas, Y.,
  Heyder, S., Peltre, O., Koep, N., Zaatiti, H., Hajri, H., Cabanes, Y.,
  Gerald, T., Chauchat, P., Shewmake, C., Brooks, D., Kainz, B., Donnat, C.,
  Holmes, S., Pennec, X.: Geomstats: {A} {Python} {Package} for {Riemannian}
  {Geometry} in {Machine} {Learning}. Journal of Machine Learning Research
  \textbf{21}(223), ~1--9 (Dec 2020)

\bibitem{misner_gravitation_1973}
Misner, C.W., Thorne, K.S., Wheeler, J.A.: Gravitation. Princeton University
  Press (1973),
  % \url{https://press.princeton.edu/books/hardcover/9780691177793/gravitation}

\bibitem{nava-yazdani_geodesic_2020}
Nava-Yazdani, E., Hege, H.C., Sullivan, T.J., von Tycowicz, C.: Geodesic
  {Analysis} in {Kendall}’s {Shape} {Space} with {Epidemiological}
  {Applications}. J Math Imaging Vis  \textbf{62}(4),  549--559 (May 2020).
  % \doi{10.1007/s10851-020-00945-w}

\bibitem{oneill_semi-riemannian_1983}
O'Neill, B.: Semi-{Riemannian} {Geometry} {With} {Applications} to
  {Relativity}. Academic Press (Jul 1983).
    %google-Books-ID: CGk1eRSjFIIC

\bibitem{pennec_parallel_2018}
Pennec, X.: Parallel {Transport} with {Pole} {Ladder}: a {Third} {Order}
  {Scheme} in {Affine} {Connection} {Spaces} which is {Exact} in {Affine}
  {Symmetric} {Spaces} (May 2018),
  % \url{https://hal.archives-ouvertes.fr/hal-01799888}, hal-01799888

\end{thebibliography}
%
 
\end{document}